\documentclass[12pt]{amsart}
\usepackage{palatino}
\usepackage{amscd,amssymb}
\usepackage{float}
\usepackage{graphicx}
\usepackage{dbnsymb}

\usepackage[cmtip,matrix,arrow]{xy}

\usepackage{graphicx}
\numberwithin{equation}{section}

\newtheorem {Proposition}[equation]     {Proposition}

\theoremstyle{definition}

\newtheorem {Question}[equation]      {Question}
\newtheorem{zremark}[equation] {Remark}
\newenvironment{Remark} {\par\footnotesize\zremark}{~\\}{\endzremark}

\newcommand{\pr} {\smallskip\noindent{\bf Proof\,\,}}

\newcommand     {\comment}[1]   {}
\newcommand     {\mute}[2] {}
\newcommand     {\printname}[1] {}

\newcommand{\labell}[1] {\label{#1}\printname{#1}}

\def    \hol    {{\operatorname{hol}}}

\def    \to     {\longrightarrow}

\def    \C      {{\mathbb C}}
\def    \R      {{\mathbb R}}

\def    \Z      {{\mathbb Z}}

\def\cH{{\bf H}}

\def    \cJ     {{\mathcal J}}
\def    \cV     {{\mathcal V}}

\def    \cA     {{\mathcal A}}

\def    \cG    {{\mathcal G}}

\def    \cN    {{\mathcal N}}

\def    \fg     {{\mathfrak g}}

\def\ba {{\bf A}}

\def    \tr     {{\rm {tr}~}}
\def    \hol   {{\rm {hol}}}
\def    \pr     {\operatorname{pr}}

\begin{document}

\title[Quantum Mechanics and Kauffman's Bracket Polynomial]{The Chern-Simons Functional Integral, Kauffman's Bracket Polynomial, and other link invariants}

\author{Jonathan Weitsman}
\thanks{Supported in part by a grant from the Simons Foundation (\# 579801)}
\address{Department of Mathematics, Northeastern University, Boston, MA 02115}
\email{j.weitsman@neu.edu}
\thanks{\today}

\begin{abstract} We study Chern-Simons Gauge Theory in axial gauge on $\R^3.$  This theory has a quadratic Lagrangian and therefore expectations can be computed nonperturbatively by explicit formulas, giving an (unbounded) linear functional on a space of polynomial functions in the gauge fields, as a mathematically well-defined avatar of the formal functional integral.  We use differential-geometric methods to extend the definition of this linear functional to expectations of products of Wilson loops corresponding to oriented links in $\R^3,$ and derive skein relations for them.  In the case $G=SU(2)$ we show that these skein relations are closely related to those of the Kauffman bracket polynomial, which is closely related to the Jones polynomial.   We also study the case of groups of higher rank.  We note that in the absence of a cubic term in the action, there is no quantization condition on the coupling $\lambda,$ which can be any complex number.  This is in line with the fact that the Jones polynomial, in contrast to the manifold invariants of Witten and Reshetikhin-Turaev, is defined for any value of the coupling.  The appearance of the parameter $e^{\frac1{2\lambda}}$ in the expectations and skein relations is also natural.  Likewise, the extension of the theory to noncompact groups presents no difficulties.
Finally we show how computations similar to ours, but for gauge fields in two dimensions, yield the Goldman bracket.

\end{abstract}

\maketitle

\tableofcontents
\section{Introduction}\labell{intro}In \cite{Witten} E. Witten showed how Chern Simons Gauge Theory is related to manifold and knot invariants, including the Jones polynomial.  Let $G$ be a compact Lie group.  Choose an irreducible representation of $G.$ This gives a representation also of the Lie algebra $\fg;$ suppose the associated trace ${\rm Tr}$ gives rise to a metric on $\fg.$\footnote{If $G$ is a classical compact group, take the defining representation.}  Let $\cA$ be the space of connections on the trivialized principal $G$ bundle on $\R^3, $ trivial outside of a compact set; then $\cA = \Omega^1_c(\R^3) \otimes \fg.$  The Chern-Simons function ${\rm CS}: \cA \to \R$ is given by 

$${\rm CS(A)} = \frac{1}{4 \pi} {\rm Tr}\int_{\R^3} A \wedge dA + \frac23 A\wedge A \wedge A.$$

We wish to study the formal functional integral

\begin{equation}\labell{pf} 
Z(\lambda) = \int_\cA DA e^{-\lambda {\rm CS}(A)}
\end{equation}

\noindent (for  $\lambda \in \C$) as well as expressions of the type
\begin{equation}\labell{exp} 
 \int_\cA DA e^{-\lambda {\rm CS}(A)} W_{C_1}(A)\dots W_{C_n}(A)\end{equation}

\noindent where $C_1,\dots C_n$ are nonintersecting oriented simple closed curves in $R^3$ and, for any such curve $C,$ the function $W_C: \cA \to \R$ is given by

$$W_C(A)= {\rm tr~} \hol_C(A)={\rm tr}_R \hol_C(A),$$

\noindent the trace of the holonomy of $A$ about $C$ in some representation $R$ of $G.$\footnote{To avoid encumbering the notation, we will choose such a representation once and for all and omit the subscript referring to it in the formulas.  In all the examples we study, the representation used for the Wilson loops will coincide with that used to define the Chern-Simons function, so ${\rm tr} = {\rm Tr}.$} Note that although the holonomy (which we denote $\hol_{C,\star}(A)$) depends on a choice of basepoint $\star \in C$, the trace does not, so we suppress the basepoint where it is irrelevant.

The formal functional integral (\ref{pf}) involves a cubic term, and therefore its mathematical interpretation is obscure.  In addition, the function ${\rm CS}$ is invariant under the action of the gauge group $\cG = {\rm Map}_c(\R^3,G)$ of maps from $\R^3$ to $G$ constant on the complement of a compact set, and acting on $\cA$ by changing the trivialization; in terms of one forms we have for $g \in \cG, A \in \cA,$

$$ g\cdot A = dg g^{-1} + g A g^{-1}.$$

It is therefore necessary to choose a section for the $\cG$ action before attempting to make sense out of (\ref{pf}) and (\ref{exp}).  In spite of these difficulties, Witten in the paper \cite{Witten} on Chern-Simons Gauge Theory studied these formal functional integrals\footnote{In the case $\lambda = i k $ where $k \in \Z.$} by physical methods and insight, and showed that the expectations (\ref{exp}) should be interpreted as link invariants giving the Jones polynomial and its generalizations, and, when $\R^3$ is replaced by a compact manifold, also predicted that the analog of (\ref{pf}) would give manifold invariants.  These were constructed topologically by Reshetikhin and Turaev \cite{rt}.  The paper \cite{Witten} also contained a number of other predictions about the invariants, all of which could be proved to be true, and was the starting point for the field of Topological Quantum Field Theory.  Nevertheless the mathematical meaning of the functional integrals remains obscure.

The purpose of this paper is to return to the functional integral (\ref{exp}) to show that in the case of $\R^3,$ the two difficulties due to gauge invariance and the appearance of a cubic term neutralize each other, and that the quadratic theory resulting from the choice of axial gauge is rich enough to contain link invariants related to the Jones polynomial and its generalizations.  

More precisely, in \cite{Witten}, the choice of section for the $\cG$ action is given by a covariant gauge of the type $d^* A = 0.$  We study instead the axial gauge condition \footnote{To obtain this gauge slice, we actually require a slightly bigger gauge group, of gauge transformations that are time independent outside a compact set and constant outside the points in $\R^3$ projecting onto some compact set in $\R^2.$  Note also that this condition fixes also the time-independent gauge transformations, due to the compact support condition on the connections.}

$$A^3 = 0.$$

In this gauge, the Faddeev-Popov ghosts decouple\footnote{See e.g. \cite{coleman}, p. 165, for the case of Yang-Mills in four dimensions; the Faddeev-Popov determinant does not depend on the Lagrangian, and is a formally infinite constant.  The calculation in three dimensions is identical.} and the cubic term in the Chern-Simons functional vanishes.  We are therefore left with a quadratic action\footnote{We absorb the factor of $2 \pi$ into the coupling.}

\begin{equation}\labell{axgaction} S(A) =  {\rm Tr}  \int_{\R^3} A^1 \frac{d}{dt} A^2\end{equation}

\noindent where $A^1, A^2$ are the components of $A$ in a fixed orthonormal basis for the planes parallel to $\R^2 \subset \R^3,$ and $t$ is the coordinate on $\R^3$ perpendicular to the plane.  We must therefore study quadratic functional integrals of the type  

\begin{equation}\labell{linka}   \int_\ba DA e^{-\lambda S(A)}W_{C_1}(A)\dots W_{C_n}(A)\end{equation}

\noindent where we denote by $\ba$ the space of connections satisfying the axial gauge condition.  Since the action $S$ is quadratic, formal expectations of polynomials in the gauge fields can be defined by
explicit formulas, and the difficulties of functional integration simplify radically.  These  quadratic
functional integrals yield expectations closely related to the link invariants as in \cite{Witten}, although some care is necessary.   First, in order to obtain the Kauffman bracket polynomial \cite{K2}\footnote{Kauffman's paper introduces a number of link invariants.  The bracket polynomial that we consider is the invariant Kauffman denotes by $\langle \cdot \rangle,$ appearing on p. 430 of \cite{K2}, and with values of the constants as on p. 432. This is related to, but not identical with, the bracket $[\cdot].$} the expectations have to be normalized in a somewhat unexpected way, with a sign inserted for the number of components and another factor for the number of crossings.  In addition, a further change has to be made, in the normalization of the expectation of a single unknotted circle, from the one expected from gauge theory.   This can be thought of in terms of a background field, and we retain a background field throughout our computations. 

Second, from a mathematical point of view, there is a distinction between the explicit formulas arising from the evaluation of expectations of polynomials in the gauge fields in the action (\ref{axgaction}) and Gaussian integrals of the type studied in Quantum Mechanics and Euclidean Quantum Field Theory \cite{gj}.  In the case of (\ref{axgaction}), the quadratic action is not positive definite, and therefore the functional integral cannot be interpreted as a continuous linear functional on the space of bounded functions on $\ba$--in other words, a measure on $\ba.$   So it is not clear that one can obtain formulas for expressions of the type (\ref{linka}) by approximation of the functions $W_C(\cdot)$ by polynomials, or exponential polynomials.  We will instead return to the heuristic functional integral and use differential-geometric methods to try to understand what expectations of the type (\ref{linka}) should look like.  The resulting formulas can be considered mathematically as definitions of the values of a linear functional associated to the formal functional integral on functions given by products of functions of the form $W_C(\cdot).$ These values are then shown to satisfy skein relations of the type studied in knot theory.  The definition of the values of the linear functional does involve returning to the formal functional integral.  However, since the formal functional integral involves a quadratic Lagrangian, the resulting formulas are exact closed form formulas, not perturbative expansions, and the skein relations have a clear mathematical interpretation as relations between values of a well-defined functional on certain functions on $\ba.$  So the main difficulties of constructive quantum field theory do not arise in this situation.\footnote{There are open questions about the nature of the linear functional, beyond the explicit formulas.  See Remark \ref{osr} and Question \ref{osq} for some of these.}

We also note that in the absence of a cubic term, there is no quantization condition on $\lambda,$ which can be any complex number.  This is in line with the fact that the Jones polynomial, in contrast to the manifold invariants of Witten and Reshetikhin-Turaev, is defined for any value of the coupling.  The appearance of the parameter $e^{\frac1{2\lambda}}$ in the expectations and skein relations is also natural.  
Likewise, the extension of the theory to noncompact groups presents no difficulties.

The fact that nonperturbative features associated with an interacting field theory may appear from the behavior of nonlinear or nonlocal functions of the fields in a free field theory is not unprecedented in quantum field theory.  One example is the vertex operator construction (see e.g. \cite{go}) which gives the expectations in the nonabelian Wess Zumino Witten model at level $1$ in terms of operators defined in a free field theory corresponding to the maximal torus of the gauge group.\footnote{In a more speculative vein, the fact that Chern Simons gauge theory is a type of string field theory leads to the question as to whether some of the features of nonperturbative string theory may arise in an analogous way. }~\footnote{The Goldman bracket \cite{goldman} is another example of the same type of phenomenon.  See Section \ref{goldmansec}.}

The structure of this paper is as follows.  In Section \ref{section2} we study Chern-Simons gauge theory on $\R^3$ in axial gauge, and study the expectations of Wilson loops by differential geometric methods applied to the heuristic functional integral.  These give relations between the expectations for Wilson loops related by certain changes in the crossings.  As a simple case we study the case $G=U(1)$ and recover the linking number, as expected in Abelian gauge theory \cite{polyakov,Witten, schwarz}.

In Section \ref{su2sec} we focus on the case $G= SU(2)$  in the defining representation.  We show how the formulas of Section \ref{section2} simplify in this case, and show how to exponentiate the action of the propagator associated with the action (\ref{axgaction}) to obtain explicit formulas for the expectations of products of Wilson loops.  In the case of the defining representation of $SU(2),$ the Wilson loops are real, which would lead us to believe that the expectations are independent of the orientations of the link components.  We show how this arises from simple formulas for traces in this representation, and as a consequence derive skein relations which after a proper normalization are precisely those of the Kauffman bracket polynomial \cite{K2}.    Iteration of those relations relates the expectation of a product of Wilson loops corresponding to any link to the expectation of a product of Wilson loops corresponding to a disjoint union of unknotted circles. Kauffman shows that his bracket polynomial is invariant under the second and third Reidemeister moves, and is closely related to the Jones polynomial, if the value of the bracket polynomial for a single unknotted circle $\langle \bigcirc \rangle$ is fixed to a particular value, which is also the value natural in conformal field theory.\footnote{To be precise, Kauffman proves that his bracket polynomial is a regular isotopy invariant of link diagrams if this choice is made.  To get an ambient isotopy invariant, an additional modification has to be made, which depends on an orientation of the link components.}  This value differs from the vacuum expectation value of the corresponding Wilson loop in gauge theory, which is $\langle \bigcirc \rangle = 2.$  The fact that axial gauge does not give quite the same answer as covariant gauge for a quantity which is in some sense a zero point energy is not entirely surprising in an anomalous theory. The difficulty can be overcome by computing expectations in the presence of a background field, and then giving a definition for the expectations of trivial Wilson loops in the background field.

In Section \ref{glnsec} we derive similar formulas for the group $GL(n,\C)$ again in the defining representation.  The group $GL(n,\C)$ is of course not a compact group, but the explicit formulas we derive in Section \ref{section2} make no use of compactness, and the formula for the action of the propagator on products of Wilson loops turns out to generalize that derived in Section \ref{su2sec} for $G=SU(2),$ with the exception of a trace term, which controls the behavior of the expectation under the first Reidemeister move.  We explain this surprise by changing the basis in which we perform the computation to a $\C$-basis for ${\mathfrak gl}(n,\C)$ which is also an $\R$-basis for ${\mathfrak u}(n),$ showing that the $GL(n,\C)$ and $U(n)$ expectations coincide.

Finally in Section \ref{goldmansec} we show how the Goldman bracket for gauge fields in two dimensions, computed in \cite{goldman} for Wilson loops restricted to the moduli space of flat connections on a Riemann surface, extends to Wilson loops considered as functions on the space of all connections. On this space, which is a symplectic manifold \cite{ab}, the Poisson bracket of Wilson loops can be computed by gauge theoretic methods similar to ours, with the symmetric propagator of Chern-Simons gauge theory in axial gauge replaced by the skew operator given by the Poisson bivector field. The resulting formulas for Goldman's Lie algebra of curves have a close resemblance to our skein relations.

\section{Chern Simons Gauge Theory on $\R^3$ in Axial Gauge}\labell{section2}

Consider then the space of connections $A \in \cA$ satisfying the condition $A^3 = 0,$ which we have denoted by $\ba.$  We wish to find concrete formulas for the expectations of products of Wilson loops in the formal functional integral given by the action $S: \ba \to \R$ given by 

 \begin{equation}\labell{axgaction2} S(A) =  {\rm Tr} \int_{\R^3} A^1 \frac{d}{dt} A^2\end{equation}

\noindent where $A^1, A^2$ are the components of $A$ in fixed coordinates in the planes parallel to $\R^2 \subset \R^3,$ $t$ is the coordinate on $\R^3$ perpendicular to the plane.  We must therefore study quadratic functional integrals of the type  

\begin{equation}\labell{linka2}  \int_\ba DA e^{-\lambda S(A)}W_{C_1}(A)\dots W_{C_n}(A)\end{equation}

The first step in carrying out these computations is to define the propagator, which is the formal inverse of the differential operator appearing in the functional $S.$  We write

\begin{equation}\labell{propagator} 
\Gamma^{ij}_{\alpha\beta} (x,t; x^\prime,t^\prime) = (g^{-1})_{\alpha\beta} \epsilon_{ij} \delta(x-x^\prime) u(t-t^\prime).
\end{equation}

Here we have chosen a basis $\{e_\alpha\}$ for $\fg,$ and define $g_{\alpha \beta} = {\rm Tr~}(e_\alpha e_\beta).$\footnote{The trace is taken in the representation used to define the Chern Simons function.  Note that if $G$ is compact we can choose a basis with $g_{\alpha\beta} = -\delta_{\alpha\beta}$. In terms of such an orthonormal basis $e_\alpha$ for $\fg,$ $S(A) = -\frac12 \sum_\alpha \int_{\R^3} A^\alpha  \wedge \frac{d}{dt} A^\alpha dt.$   (The factor of $\frac12$ arises from the sum over cotangent indices replacing the choice $A_1 \frac{dA_2}{dt}.$)}.  Also, $x,x^\prime$ are the coordinates on $\R^2$ (and their pullbacks by projection to $\R^3$), $t,t^\prime \in \R,$ and\footnote{The addition of a constant to $u,$ although in principle allowable, results in no change to any expectation, due to the fact that $\epsilon_{ij}$ is skew.} 

$$u(s) =  \left\{
	\begin{array}{ll}
		\frac12  & \mbox{if } s \geq 0 \\
		-\frac12 & \mbox{if } s< 0
	\end{array}
\right.$$

We now define some functions on $\ba. $  Given $A \in \ba,$ $f \in C_c^\infty(\R^3) \otimes \R^2 \otimes \fg,$ write, using an orthonormal basis for $\fg,$
 
$$A(f) = \int_{\R^3} \sum_\alpha \sum_i A_\alpha^i(x,t) f_\alpha^i(x,t) d^2x dt.$$  This expression gives a function on $\ba$ we continue to denote by $A(f).$

Then by analogy to the usual Gaussian integrals, in either finite or infinite dimensions, we use the propagator (\ref{propagator}) to define a linear functional on quadratic polynomials in such functions $A(f), A(\tilde{f})$ for $f,\tilde{f}  \in C_c^\infty(\R^3) \otimes \R^2 \otimes \fg,$ by 

$$\langle A(f) A(\tilde{f}) \rangle_0 = \lambda^{-1} \sum_{\alpha,\beta} \sum_{i,j} \int_{\R^3} (g^{-1})_{\alpha\beta}\epsilon_{ij} f_\alpha^i(x,t) \tilde{f}_\beta^j(x,t^\prime) u(t-t^\prime) d^2x dt.$$

We then extend the definition of the linear functional $\langle\cdot\rangle_0$ to all polynomials by analogy with 
Wick's theorem.  Given  $f_1,\dots,f_n  \in C_c^\infty(\R^3) \otimes \R^2 \otimes \fg,$ where $n$ is even, let

\begin{equation}\labell{wick}
\langle A(f_1) \dots A(f_n) \rangle_0 = \frac{1}{2^{\frac{n}{2}}(\frac{n}{2}!)}\sum_{\sigma \in \Sigma_n} \langle A(f_{\sigma(1)} )  A(f_{\sigma(2)} )\rangle_0\dots \langle A(f_{\sigma(n-1)} )  A(f_{\sigma(n)})\rangle_0 \end{equation}

\noindent where the sum is taken over all elements $\sigma$ in the permutation group $\Sigma_n.$  The functional $\langle\cdot\rangle_0$ is defined to be zero on all odd degree monomials, and we set

$$\langle 1 \rangle_0 = 1.$$

Another way of expressing these formulas is by writing

\begin{equation}\labell{pertser}\langle P \rangle_0 = e^{\frac1{2\lambda} \sum_{\alpha,\beta} \sum_{i,j}   \int_{\R^3 \times \R^3} dx dy \frac{\delta}{\delta A_\alpha^i(x)} \Gamma_{\alpha \beta}^{ij}(x,y)  \frac{\delta}{\delta A_\beta^j(y)}} P|_{A=0},\end{equation}

\noindent where $P$ is any polynomial given as a sum of monomials of the form $A(f_1) \dots A(f_n)$ as above, and the operator $\Gamma$ is given in (\ref{propagator}).\footnote{Note that evaluation of the linear functional $\langle \cdot \rangle,$ in analogy to the case of Gaussian integrals, is given by formula (\ref{pertser}), which involves only functional {\em derivatives}.  Differentiation in function spaces, unlike integration on function spaces, is straightforward to define.}

\begin{Remark}\labell{osr}{\bf Some comments on Reflection Positivity.} 

Let $G$ be a compact Lie group and suppose the trace ${\rm Tr}$ allows us to choose an orthonormal basis $\{e_\alpha\}$ for $\fg$ so that ${\rm Tr} (e_\alpha e_\beta) = - \delta_{\alpha\beta}.$\footnote{This is the case for a compact simple Lie group with the trace in the adjoint representation.}   Suppose that ${\rm Supp}(f) \subset \R^3_+$ and ${\rm Supp}(g) \subset \R^3_-,$ where $\R^3_\pm = \{(x,y,z) \in \R^3:  \pm z > 0\}.$  Then\footnote{The factor of $\frac12$ in the following equation arises from the function $u.$ The sign arises from the difference between the metric and the trace on the Lie algebra $\fg.$}

\begin{equation}\labell{os1}
\langle A(f) A(g) \rangle_0 = -  \frac1{2\lambda} \langle \pi(f) , \pi(rg)\rangle_{L_2(\R^2) \otimes \R^2 \otimes \fg}\end{equation}

\noindent where $\pi: C_c^\infty(\R^3) \otimes \R^2 \otimes \fg \to C_c^\infty(\R^2) \otimes \R^2 \otimes \fg$ is integration along the $t$-direction\footnote{This integration along lines is part of the Radon transform.  See e.g. \cite{helg2}.}  and $r$ is the rotation by $\pi/2$ in the cotangent directions in the planes parallel to $\R^2$:  explicitly

$$(\pi(f))(x) = \int_\R f(x,t) dt$$

\noindent and

$$\begin{array}{ll}
		(r f)^2_\alpha  & =  - f^1_\alpha \\
		(rf)^1_\alpha  & = f^2_\alpha
	\end{array}$$

\noindent for $f \in C_c^\infty(\R^3) \otimes \R^2 \otimes \fg.$ 

Now consider the case where ${\rm Supp}(f) \subset \R^3_+$ as above, and define the reflection of $f$ by 

$$Rf(x,t) = f(x, -t),$$

\noindent so that ${\rm Supp}(Rf) \subset \R^3_-.$  Then for $\lambda > 0,$ the expression 

$$\langle A(f) A(rRf) \rangle_0$$

\noindent defines a positive semidefinite inner product on $\cV_1 = C_c^\infty(\R^3_+) \otimes \R^2 \otimes \fg.$  As in \cite{os,gj}, let $\cN_1$ be the set of null elements of $\cV_1$ for this inner product, and define the one particle Hilbert space $\cH_1$ as the completion of the quotient $\cV_1/\cN_1$ in the resulting positive definite inner product.  Then

\begin{Proposition}\labell{ost}

For $\lambda > 0,$ the linear functional $\langle\cdot\rangle_0$ gives a reflection positive functional on $\cV_1=C_c^\infty(\R^3_+) \otimes \R^2 \otimes \fg;$ we have

\begin{equation}\labell{os2}
\langle A(f) A(rRf) \rangle_0 =  \frac1{2\lambda} ||\pi(f)||^2_{L_2(\R^2) \otimes \R^2 \otimes \fg}\end{equation}

\noindent for $f \in \cV_1.$ The one particle Hilbert space is given by\footnote{The $L_2$ metric is scaled by $ 2 \lambda.$}

$$\cH_1 = {L_2(\R^2) \otimes \R^2 \otimes \fg}.$$

\end{Proposition} 

However, unlike in the case of free Bosonic fields, the propagator $\Gamma$ is not invariant under the action of the element $Rr \in O(3),$\footnote{One manifestation of that is that if $f,g \in \cV_1,$ then $\langle A(f) A(g) \rangle_0 = - \langle A(rRf) A(rRg)\rangle_0,$ so that Wick ordering is not reflection invariant.  This is due to the fact that the action (\ref{axgaction2}) is not invariant under orientation reversing diffeomorphisms.} so that reflection positivity for linear functions does not guarantee reflection positivity for polynomials in these linear functions.  Thus our construction of the one particle subspace $\cH_1$ does not extend in a straightforward way to the construction of a Hilbert space for the quantum field theory.  Difficulties of this type are known to arise in theories with one derivative in the free field action, and various methods such as doubling the fields \cite{osf} or restricting to a subalgebra \cite{jp} have been used to develop useful forms of reflection positivity in such cases.

\begin{Question}\labell{osq} Is there a form of reflection positivity for Chern-Simons quantum field theory in axial gauge?\end{Question} \end{Remark}

\subsection{The expectation of a product of Wilson loops} Our next task is to try to make sense of the value of the functional $\langle \cdot \rangle_0$ on functions of the form 
$W_{C_1}\dots W_{C_n}.$  This is not entirely straightforward, since the functions $W_{C}$ are not in the domain of definition of the functional $\langle \cdot \rangle_0;$ and since this functional is not known to be bounded, it is not clear that any approximation scheme of the $W_C$ by polynomials will be useful. Instead, we return to the formal functional integral (\ref{exp}) to find an independent definition for the expectations of such objects.

\begin{Remark} To the extent we expect a linear functional on the algebra of bounded functions on $\ba,$ such as $\langle \cdot \rangle_0$ (or a linear functional on the bounded functions on $\cA,$ for that matter) to give the Jones polynomial, we may expect that unboundedness of the functional is essential.  Take $G = SU(2)$ in the defining representation, and a single Wilson loop $W_C;$ $W_C$ is a bounded function on the space of connections $\ba,$ and $||W_C||_{\infty}=2.$  However, one can find knots with arbitrarily large Jones polynomial.  (One example is the connected sum of multiple copies of the trefoil.)\footnote{I would like to thank Fabian Ruehle for this example.}  So it is not possible for the Jones polynomial to arise as the expectation of $W_C$ in any measure on $\ba, $ or $\cA.$\footnote{Note that in quantum field theory, the measures corresponding to the Euclidean action tend to be defined on the analogs of completions of $\cA$ to spaces of distributions.  The definition of functions of the type $W_C$ on such distributional connections is not entirely straightforward, and may require renormalization.  This may be another way of looking at the same issue.}  \end{Remark}

We begin by studying connections on the circle.  Let $A$ be a connection on the trivial principal $G-$bundle on $C=S^1,$ with its standard orientation, and let $v \in \Omega^1(S^1)\otimes \fg.$  We compute

\begin{equation}\labell{conn1}\frac{d}{ds}\vert_{s=0} {\rm tr~} \hol_C (A + s v) = \int_0^{2 \pi} dt~ {\rm tr~}( \hol_{C,t} (A) v(t) ).\end{equation}

Here we have written $\hol_{C,t} (A) \in G $ for the holonomy of the connection $A$ along the circle from the point $t$ to itself in the direction given by the standard orientation, and we consider both $G$ and $\fg$ as subsets of some space $M_n(\C)$ of $n \times n$ matrices corresponding to the chosen representation.

In terms of functional derivatives, this equation (\ref{conn1}) can be written

\begin{equation}\labell{conn}
\frac{\delta}{\delta A_\alpha(t)} {\rm tr~} \hol_C(A) = {\rm tr~} (p( \hol_{C,t} (A)) e_\alpha)\end{equation}

\noindent where $p: G \to \fg$ is the projection arising from the representation we chose for $G$ (inside some $GL(n,\C)$) and the inner product arising from the corresponding trace.

For future use we also record this functional derivative for the case of a connection on $[0,1].$  We have (using again the representation in $GL(n,\C)$)\footnote{The factor of $\frac12$ in the equations below is due to the derivative being taken at the endpoint, and can be seen by applying the analog of (\ref{conn1}) to a series of approximations of the delta function by smooth functions.  If the connection on a circle is thought of as being a connection on the interval with endpoints identified, the factors of $\frac12$ arising from the endpoints add to give (\ref{conn1}).}

\begin{equation}\labell{connpath}
\frac{\delta}{\delta A_\alpha(0)}   \hol_{[0,1]} (A) = \frac12 e_\alpha  \hol_{[0,1]} (A)   \in  M_n(\C). \end{equation}

\noindent  Likewise
 
\begin{equation}\labell{connpathe}
\frac{\delta}{\delta A_\alpha(1)}   \hol_{[0,1]} (A) =  \frac12 \hol_{[0,1]} (A) e_\alpha  \in M_n(\C) \end{equation}

Given an embedded path, or curve, in a manifold, and a connection on the trivialized principal $G$ bundle on the manifold, these formulas give expressions for the derivatives in the tangent direction to the curve of the holonomy of the pulled back connection.

We now apply the formula (\ref{conn}) to calculate the expectation of a function given by a product of Wilson loops.  Suppose we are given a collection $C_1,\dots,C_n$ of nonintersecting oriented simple closed curves in $\R^3.$  Suppose that the projections of these curves onto the plane $\R^2$ is a family of curves intersecting transversally 
at a finite number of double points $x_1,\dots,x_k.$  Denote the inverse images
of the double points in $\R^3$ by $(x_i,t_i^\pm).$ Then the formal functional integral (\ref{linka}) leads us to consider, in analogy to (\ref{pertser}),  the expression

\begin{equation}\labell{exp1}
\langle W_{C_1} \dots W_{C_n}\rangle = \prod_{i=1}^k e^{\frac{1}{2\lambda} D_i } W_{C_1}(A) \dots W_{C_n}(A) \end{equation}

\noindent where  
 
\begin{equation}\labell{op} D_i=   \sum_{\alpha,\beta} \sum_{a,b} \sum_{s,s^\prime \in \{t_i^+,t_i^-\}}(g^{-1})_{\alpha\beta} \epsilon_{ab}  u(s - s^\prime)   \frac{\delta}{\delta A_\alpha^a(x_i,s)}\frac{\delta}{\delta A_\beta^b(x_i,s^\prime)}. \end{equation}

\noindent The formal Gaussian integral (\ref{linka}) would lead to the evaluation of $\langle W_{C_1} \dots W_{C_n}\rangle $ at $A=0;$ we will see it is useful to keep track of $\langle W_{C_1} \dots W_{C_n}\rangle$ as a function of $A,$ which may be thought of as a background gauge field.  

\begin{Remark} We emphasize that by the expression (\ref{exp1}), and all the explicit formulas that follow, we mean the repeated application of the operators $D_i,$ and their exponentials, to products of the $W_{C_i},$ as functions of the background gauge field $A,$ not to the evaluation of those expressions at $A=0.$ In fact, it will be convenient to evaluate such functions related by a skein relation on connections which are pullbacks from the plane, {\em after} application of the operators in (\ref{exp1}). \end{Remark}

As a first step, we  compute $$D_i (W_{C_1} \dots W_{C_n}).$$

Since the operator $D_i$ is second order, it suffices to consider its action on linear and quadratic expressions in the functions  $W_{C_i}.$

Consider first the action of one operator $D_i$ on one Wilson loop $W_C.$  This corresponds to one pair of points $(x,t^\pm) = (x_i,t^\pm_i)$ corresponding to one of the double points $x = x_i$ of the projection of $C$ to the plane. 
Removing the two points $(x,t^\pm)$ cuts $C$ into two oriented paths with closures $X,Y.$\footnote{Of course the paths $X$ and $Y$ may contain other points projecting to intersection points.}
 Choose a basis for the tangent space to $\R^2$ at $x$ given by the tangent vectors to the projections of the oriented paths $X$ and $Y$ at their starting points at $x,$ and orient the basis using the usual orientation of the plane; denote the oriented basis by $(X_1,X_2).$  By exchanging the names of $X$ and $Y$ if necessary, we may assume $X_1$ is the tangent vector to $X$ at its initial point and $X_2$ is the tangent to $Y$ at its initial point.  We also rename $t^\pm$ so that $X$ travels from $(x,t^+)$ to $(x,t^-)$ and $Y$ travels from $(x,t^-)$ to $(x,t^+).$\footnote{ Note that in our convention $t^+ > t^-$ for an overcrossing and $t^- > t^+$ for an undercrossing. I do apologize.  The alternative would be to orient the tangent vectors to $X$ and $Y$ with a right handed orientation for an overcrossing and a left handed orientation for an undercrossing.} (See Figure 1.)  
 
 \begin{figure}\label{figure1}
       \includegraphics[width=0.65\textwidth]{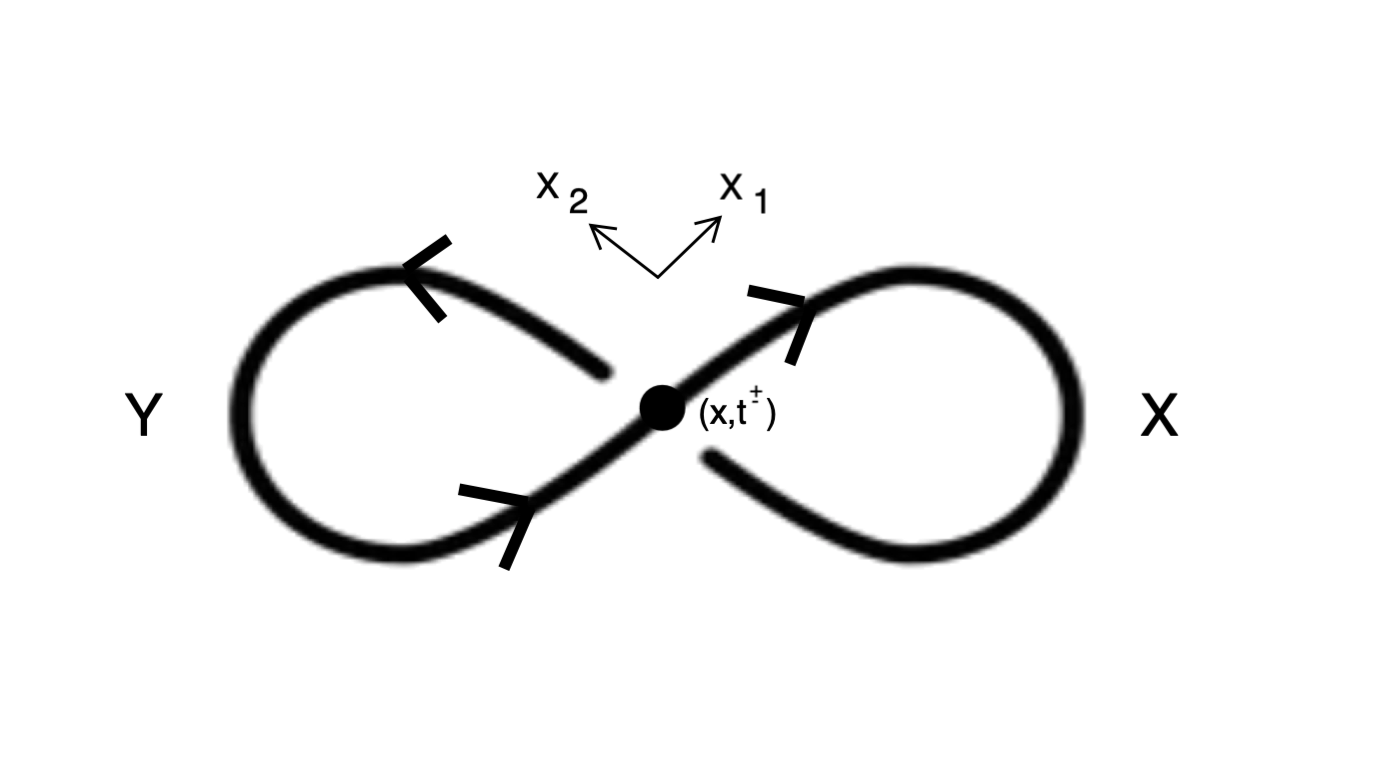}
    \caption{The curve $C$ cut into two paths $X,Y$ by removing $(x,t^\pm)$}
   \end{figure}
 
 We wish to compute\footnote{The factor of 4 is due to the choice of $a,b$ and $s,s^\prime$ in the operator (\ref{op}) being fixed.} 

$$  4\sum_{\alpha,\beta} (g^{-1})_{\alpha\beta} u(t^+ - t^-) \frac{\delta}{\delta A_\alpha^1(x,t^+)}\frac{\delta}{\delta A_\beta^2(x ,t ^-)} W_C(A)$$

\noindent where the superscripts $1,2$ refer to the oriented basis $(X_1,X_2)$ for the tangent space at $x.$\footnote{Strictly speaking the propagator (\ref{propagator}) involved the cotangent indices in the directions given by the coordinate axes on $\R^2,$ but the final expressions are independent of the coordinates chosen, and depend only on the orientation on $\R^2,$ so that we may write a similar formula in any coordinate system.  In particular, there is no need for the two branches of $C$ above $x$ to meet at right angles.}

Applying formulas
(\ref{connpath})-(\ref{connpathe}), we have

\begin{equation}\labell{linear}
 4 \sum_{\alpha,\beta} (g^{-1})_{\alpha\beta}     u(t^+ - t^-) \frac{\delta}{\delta A_\alpha^1(x,t^+)}\frac{\delta}{\delta A_\beta^2(x ,t ^-)}  W_C(A) = 4 u(t^+ - t^-) \sum_{\alpha,\beta} (g^{-1})_{\alpha\beta}  {\rm tr~} e_\alpha\hol_X(A) e_\beta \hol_Y(A) 
\end{equation}

\noindent where $\hol_X(A)$ and $\hol_Y(A)$ denote the holonomy of $A$ along the oriented paths $X,Y.$\footnote{The two factors of $\frac12$ appearing in 
(\ref{connpath})-(\ref{connpathe}) are cancelled by the fact that each derivative can be placed at either endpoint of the paths $X$ and $Y$, giving two factors of 2.}

Similarly, for the quadratic case, we compute
 
$$ 4 \sum_{\alpha,\beta} (g^{-1})_{\alpha\beta}   u(t^+ - t^-) \frac{\delta}{\delta A_\alpha^1(x,t^+)}\frac{\delta}{\delta A_\beta^2(x ,t ^-)} W_{C}(A) W_{C^\prime}(A)$$

\noindent where $(x,t^\pm)$ are points of the curves $C,C^\prime$ respectively.\footnote{Again, as in the case of a single curve, the projections of the curves $C,C'$ to the plane may intersect at other points as well.}  We obtain 
 
\begin{equation}\begin{aligned}\labell{quad}
4\sum_{\alpha,\beta} (g^{-1})_{\alpha\beta}   u(t^+ - t^-) \frac{\delta}{\delta A_\alpha^1(x,t^+)}\frac{\delta}{\delta A_\beta^2(x ,t ^-)}  W_{C}(A) W_{C^\prime}(A)=&\\
4\sum_{\alpha,\beta} (g^{-1})_{\alpha\beta}  u(t^+ - t^-)({\rm tr~} e_\alpha \hol_{C,(x,t^+)}(A))( {\rm tr~} e_\beta  \hol_{C^\prime,(x,t^-)}(A) ). \end{aligned}\end{equation}
 \noindent (See Figure 2.)

 \begin{figure}\label{figure2}
       \includegraphics[width=0.5\textwidth]{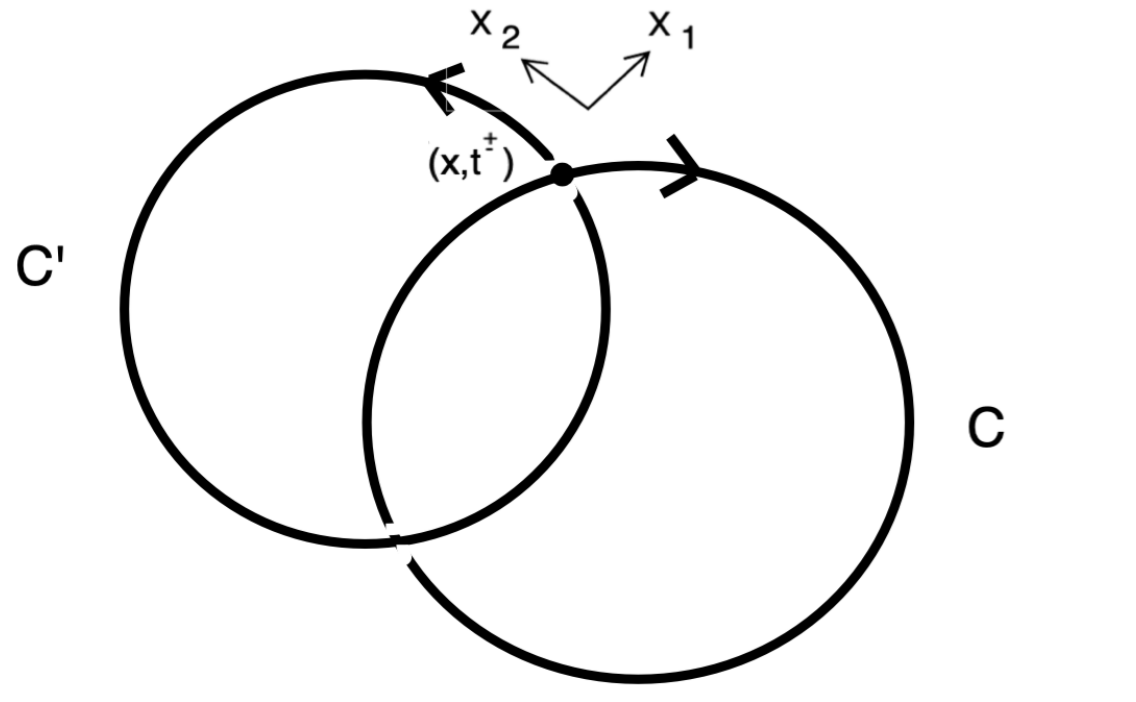}
    \caption{The curves $C$ and $C^\prime$ and the points $(x,t^\pm)$}
   \end{figure} 
   
In the next two sections, we show how to develop effective formulas for the expressions appearing on the right hand sides of equations (\ref{linear}) and (\ref{quad}).\footnote{The formulas (\ref{linear}) and (\ref{quad}) appear in some way to represent a theory occurring on a space $\R^2 \times \{0,1\}$ consisting two disjoint copies of the plane, with interactions occurring at pairs of points $(x,0)$ and $(x,1)$ with $x \in \R^2.$  This is in some ways reminiscent of Connes' noncommutative geometry models of gauge theory Lagrangians.  See \cite{Connes}, pp. 561 ff.}  We begin with the simplest case, $G = U(1).$

 \subsection{The $U(1)$ case}\labell{u1sec}    Take $G=U(1)$ in the standard one dimensional representation.  In this case there is one basis element $e_1 = i$ for the Lie algebra $i\R,$ $g_{11} = -1,$ and equation (\ref{linear}) simplifies to
 
\begin{equation}\labell{linearu1}
-4 u(t^+ - t^-) \frac{\delta}{\delta A^1(x,t^+)}\frac{\delta}{\delta A^2(x ,t ^-)}  W_C(A) = -  4 u(t^+ - t^-) i  \hol_X(A) i \hol_Y(A) =2 \epsilon W_C(A)
\end{equation}

\noindent where $\epsilon = +1$ for an overcrossing and $\epsilon = -1$ for an undercrossing; so that $u(t^+ - t^-) = \frac12 \epsilon.$  Similarly equation (\ref{quad}) simplifies to 
 
\begin{equation}\labell{quadu1}
-4 u(t^+ - t^-) \frac{\delta}{\delta A^1(x,t^+)}\frac{\delta}{\delta A ^2(x ,t ^-)}  W_{C}(A) W_{C^\prime}(A)= 
-  2\epsilon i  \hol_{C,(x,t^+)}(A) i \hol_{C^\prime,(x,t^-)}(A) =  2\epsilon W_{C}(A) W_{C^\prime}(A). \end{equation}

\noindent with the same definition for $\epsilon.$  

Note that in both cases, the effect of the application of the operator $D$ at a crossing, whether of a strand with itself or between two strands, is the same; a similar phenomenon will occur for other groups as well, and is the reason that we obtain skein relations.  In any event, the effect of the application of $D_i$ is simply multiplication by $\pm 2,$ so that

\begin{equation} \prod_{i=1}^k e^{\frac{1}{2\lambda} D_i } W_{C_1}(A) \dots W_{C_n}(A) =  
\prod_{i=1}^k e^{\frac{1}{\lambda} \epsilon_i} W_{C_1}(A) \dots W_{C_n}(A)=
 e^{\frac{1}{\lambda} w_{C_1,\dots,C_n}} W_{C_1}(A) \dots W_{C_n}(A)
\end{equation}

\noindent where $\epsilon_i = +1$ if the crossing at $(x_i,t_i^\pm)$ is an overcrossing and 
$\epsilon_i = -1$ if the crossing at $(x_i,t_i^\pm)$ is an undercrossing; and we define the writhe $w_{C_1,\dots,C_n}$ of the oriented link consisting of ${C_1,\dots,C_n}$ by 

$$w_{C_1,\dots,C_n} = \sum_i \epsilon_i.$$

Thus our computation gives the linking number of the link given by ${C_1,\dots,C_n}.$  This agrees with computations in Abelian gauge theory in covariant gauge \cite{polyakov,Witten, schwarz}.

\section{$G = SU(2)$}\labell{su2sec}

Now let $G=SU(2)$ and choose the defining (two-dimensional) representation for both the Chern-Simons term and the Wilson loops.  Our goal is to obtain an explicit expression for (\ref{exp1}) in this case.  We start by computing the effect of the action of a single operator $D_i$ at a single double point.

Consider first the linear case, where the points $(x_i,t_i^\pm)$ occur on a single curve $C.$  We write $x = x_i, t^\pm = t_i^\pm.$  Then, as in equation (\ref{linear}), we have 

$$ {\rm tr~} \hol_C(A) = {\rm tr~} U^{12} V^{21}$$

\noindent where in the notation of (\ref{linear})

$$U^{12} = {\rm tr~} \hol_X(A)$$
$$V^{21} = {\rm tr~} \hol_Y(A)$$

\noindent where $X,Y$ are the oriented paths obtained by deleting the two points $(x,t^\pm)$ from $C.$   (The notation $U^{12}, V^{21}$ for the holonomies indicates the tangent directions of the curves $X,Y,$ at their initial and final points.  See Figure 1.)

Then applying the operator $D_i= D$, and choosing an orthonormal basis $\{e_\alpha\}$ for $\fg,$\footnote{In the metric arising from the trace.} so that $g_{\alpha\beta} = - \delta_{\alpha\beta},$ we have

$$D W_C(A) = - 2 \epsilon \sum_\alpha {\rm tr~} U^{12} e_\alpha V^{21} e_\alpha$$
\mute{factor of 4 from equation, divide by the 2 in the def of $u(t)$.}
\noindent where $\epsilon = \pm 1$ as before depending on whether the double point corresponds to an overcrossing or an undercrossing.  

Let us choose the orthonormal basis $e_\alpha = \frac{1}{\sqrt{2}} \sigma_\alpha, \alpha = 1, 2, 3,$ for ${\mathfrak su}(2),$ where $\sigma_\alpha, \alpha = 1, 2, 3$  are the Pauli matrices, which satisfy the relations $\sigma_\alpha^2 = -1$ and $\sigma_\alpha \sigma_\beta = \sum_\gamma \epsilon_{\alpha\beta\gamma} \sigma_\gamma$.  We also choose coordinates on $\R^2$ near the point $x$ so that the coordinate axes lie along the directions taken by the projection of the curve $C$ to $\R^2.$  We obtain

$$D W_C(A) = - \epsilon \sum_\alpha {\rm tr~} U^{12} \sigma_\alpha V^{21} \sigma_\alpha.$$

Now if $g \in SU(2),$ we have

\begin{equation}\labell{su2factt}
\sum_\alpha \sigma_\alpha g \sigma_\alpha = g - 2I{\rm tr~}g 
\end{equation}

\noindent where $I$ is the identity matrix.  So

\begin{equation}\labell{linearsu2}
D W_C(A) = \epsilon(- {\rm tr~} U^{12} V^{21} + 2 {\rm tr~} U^{12} {\rm tr~} V^{21}).
\end{equation}

Note that, intuitively, the terms appearing on the right hand side of (\ref{linearsu2}) correspond to curves formed from the paths $X$ and $Y;$ in the first term we glue $X$ and $Y$ back to form $C$ while in the second term we take the holonomy about the closed paths formed from $X$ and $Y$ by "jumping" from $(x,t^+)$ to $(x,t^-).$  We will see that this is a recurring phenomenon responsible for the skein relations.

We perform a similar computation in the quadratic case, where we have a single double point arising from two points $(x,t^\pm)$ occurring on two curves $C,C^\prime.$ We again choose coordinates on $\R^2$ near the point $x$ where the coordinate axes lie along the directions taken by the projection of the curve $C,C^\prime$ to $\R^2.$  We then have

\begin{equation}\labell{quadsu2}D W_C(A) W_{C^\prime}(A) = - \epsilon  \sum_\alpha{\rm tr~} U^{11} \sigma_\alpha {\rm tr~} V^{22} \sigma_\alpha = \epsilon(-{\rm tr~} U^{11} {\rm tr~} V^{22} + 2 {\rm tr~} U^{11} V^{22}), \end{equation}

\noindent where we have written $U^{11}= \hol_{C,(x,t^+)} (A)$ and $V^{22}= \hol_{C,(x,t^-)} (A);$ the notation $U^{11}, V^{22}$ for the holonomies again indicates the tangent directions of the curves $C,C^\prime,$ at their initial and final points, both of which project down to the double point $x.$ See Figure 2.

Note that to obtain the second equality in (\ref{quadsu2}) we have used the relation

$$-\sum_\alpha {\rm tr~} U \sigma_\alpha {\rm tr~} V \sigma_\alpha   = 2 {\rm tr~} UV - {\rm tr~}U {\rm tr~} V$$

\noindent for $U,V \in SU(2).$

Note also that the terms appearing on the right hand side of the equation (\ref{quadsu2}) arise from the curves $C, C^\prime$ either by leaving them alone or by gluing $C$ to $C^\prime$ at the points $(x,t^\pm)$ by a jump from $(x,t^-)$ to $(x,t^+),$ just as in (\ref{linearsu2}).

Therefore in either the linear or the quadratic case, we obtain the relations

\begin{equation}\labell{urskeino}
D ( \overcrossing)= -1 ( \overcrossing)+ 2 ( \upupsmoothing)\end{equation}

\begin{equation}\labell{urskeinu}
D (\undercrossing)= 1 (\undercrossing)- 2 (\upupsmoothing).\end{equation}

\noindent where we have used the standard notation in knot theory for relations among links involving one crossing at a time.  Note that in our case the relations are among functions on the space of connections, not only their values at the trivial connection.

Note that the $\upupsmoothing$ term in (\ref{urskeino}) and (\ref{urskeinu}) actually corresponds to holonomies of two curves that have jumps at $(x,t^\pm),$ not to curves that are separated by a gap; this will be important when we compute the action of higher powers of the operators $D_i.$  The precise formulas are equations (\ref{linearsu2}) and (\ref{quadsu2}).\footnote{The difference between the smoothing of an undercrossing and that of an overcrossing can be dealt with in the final formulas (such as (\ref{ku1}) and (\ref{kbracket})) by evaluating the resulting functions on both sides of any skein relation, obtained after computation of the action of all the operators $D_i,$ at connections which are pullbacks from the plane.  More geometrically, think of the closed loops formed from the paths corresponding to $U^{12}, V^{21}$ and the like by connecting the initial and final points by a nearly vertical strand on which we define the associated holonomy to be the identity.}

A similar computation shows that, in the linear case,

$$D {\rm tr~} \hol_X(A) {\rm tr~} \hol_Y(A) = \epsilon( {\rm tr~} \hol_X(A) {\rm tr~} \hol_Y(A)+ W_C(A))$$

\noindent or

$$D \tr U^{12} \tr V^{21} = \epsilon  (\tr U^{12} \tr V^{21} +  \tr U^{12} V^{21}),$$

\noindent while in the quadratic case

$$D {\rm tr~}\hol_{C,(x,t^-)} (A) \hol_{C^\prime,(x,t^+)}(A) = \epsilon( {\rm tr~}\hol_{C,(x,t^-)} (A) \hol_{C^\prime,(x,t^+)}(A) + W_C(A) W_{C^\prime}(A)) $$

\noindent or

$$D \tr U^{11}V^{22} = \epsilon  (\tr U^{11} V^{22} +  \tr U^{11} \tr V^{22}).$$

These results show that for either the linear or quadratic case, the effect of the application of the operator $D$ is given by

\begin{equation}
 D
  \begin{bmatrix}
   \overcrossing  \\
 \upupsmoothing 
  \end{bmatrix}
 =
M
   \begin{bmatrix}
     \overcrossing    \\
   \upupsmoothing 
  \end{bmatrix}
\end{equation}

\begin{equation}
 D
  \begin{bmatrix}
     \undercrossing    \\
   \upupsmoothing 
    \end{bmatrix}
 =
-M
   \begin{bmatrix}
   \undercrossing   \\
  \upupsmoothing 
  \end{bmatrix}
\end{equation}

\noindent where $M$ is the $2 \times 2$ matrix 

\begin{equation}\labell{m}
M = \begin{bmatrix}
    -1 & 1   \\
   2 & 1
  \end{bmatrix}\end{equation}
  
To compute expressions of the type appearing in equation (\ref{exp1}), we must exponentiate the matrix $M.$  Since $M$ is traceless, $M^2$ is a multiple of the identity matrix; in fact $M^2 = 3 I,$ and we have

\begin{equation}\labell{expm}
e^{\beta M} = \cosh(\sqrt{3} \beta) I + \frac{1}{\sqrt{3}}\sinh(\sqrt{3} \beta) M 
\end{equation}

\noindent where we have written $\beta = \frac{1}{2\lambda}.$  By computing $$e^{\beta M}  \begin{bmatrix}
   \overcrossing   \\
   0
  \end{bmatrix}$$ we therefore obtain the skein relation

\begin{equation}\labell{skeinsu2over}
\langle \overcrossing \rangle = \Big(\cosh(\sqrt{3} \beta) - \frac{1}{\sqrt{3}}\sinh(\sqrt{3} \beta)\Big) (\doublepoint )
+ \Big(\frac{2}{\sqrt{3}}\sinh(\sqrt{3} \beta)\Big) (\upupsmoothing)\end{equation}

and similarly

\begin{equation}\labell{skeinsu2under}
\langle \undercrossing \rangle =  \Big(\cosh(\sqrt{3} \beta)+\frac{1}{\sqrt{3}}\sinh(\sqrt{3} \beta)\Big) (\doublepoint )
- \Big(\frac{2}{\sqrt{3}}\sinh(\sqrt{3} \beta)\Big) (\upupsmoothing).\end{equation}

Here we have used the notation ($\doublepoint$) for a "bare" crossing; one where all the appropriate operators $D_i$ have already been applied.  This notation has been used by Kauffman \cite{kauffman} to denote a virtual crossing.  We see that our quantum field theory automatically includes Kauffman's virtual knots.  Similarly ($\upupsmoothing$) for the "honest" smoothed intersection (in our case involving a "jump").

Adding (\ref{skeinsu2over}) and (\ref{skeinsu2under}) we obtain

\begin{equation} \labell{su2+}
\langle \overcrossing \rangle + \langle \undercrossing\rangle = \Big( 2 \cosh(\sqrt{3} \beta) \Big)(\doublepoint)\end{equation}

\noindent while subtracting, we obtain

\begin{equation} \labell{su2-}
\langle \overcrossing \rangle - \langle \undercrossing\rangle = \Big( \frac2{\sqrt{3}} \sinh(\sqrt{3} \beta) \Big)\Big(- (\doublepoint) + 2 (\upupsmoothing)\Big)\end{equation}

It is possible to use the  expressions (\ref{su2+}) and (\ref{su2-}) to eliminate the virtual crossings and obtain a skein relation involving only overcrossings, undercrossings, and smoothings.  We will do this in Section \ref{glnsec} where we consider the groups $GL(n)$ and $U(n).$  

In the $SU(2)$ case however, there are further simplifications due to the fact that the trace of a matrix in $SU(2)$ is real; for $g \in SU(2),$

$$ \tr g = \tr g^{-1},$$

\noindent so that the expectations are independent of orientation, and therefore are associated to unoriented links.\footnote{The resulting expectations will turn out to be regular isotopy invariants of link diagrams.  They are related to ambient isotopy invariants of knots by a twist that is defined using a choice of orientation of the link components.  See \cite{K2} or further discussion below.}

In particular, we have the formula

\begin{equation}\labell{su2fact}g + g^{-1} =  (\tr g)I.\end{equation}

We will now show how this allows us to convert the relations (\ref{skeinsu2over}) and (\ref{skeinsu2under}) into relations among unoriented link diagrams. 

To do this we apply (\ref{su2fact}) to the expressions ${\rm tr~} U^{12} V^{21} $ and ${\rm tr~} U^{11} \tr V^{22}$ corresponding to the terms involving a crossing $\langle \doublepoint\rangle $ in   (\ref{skeinsu2over}) and (\ref{skeinsu2under}).  We have

\begin{equation}\labell{su2fact1}\tr U^{12} V^{21} = \tr\Big(  (\tr U^{12}) -(U^{12})^{-1} \Big) V^{21} = \tr U^{12} \tr V^{21} - \tr (U^{12})^{-1} V^{21} \end{equation}

\noindent while

\begin{equation}\labell{su2fact2} \tr U^{11} \tr V^{22}  = \tr \Big( (\tr U^{11}) V^{22} \Big)
=\tr \Big( (U^{11} + (U^{11})^{-1})V^{22} \Big) = 
\tr U^{11}  V^{22} +  \tr (U^{11})^{-1} V^{22}.\end{equation}

To remove the sign discrepancy between the right hand sides of (\ref{su2fact1}) and (\ref{su2fact2}), we adjust the normalization of our expectations by a factor of $(-1)^{j}$ where $j$ is the number of components of the resulting link.  We then obtain, for either linear or quadratic functions, the relation
\footnote{To avoid encumbering the formulas, we abuse notation by using the same notation for the expectation with the modified normalization.}

\begin{equation}\labell{su2miracle}
(\doublepoint) = -\Big( (\smoothing) + (\hsmoothing) \Big)
\end{equation}

\noindent so that the relations (\ref{skeinsu2over}) and (\ref{skeinsu2under})  become\footnote{  Note that the change of normalization inserts a sign in the second (smoothing) terms of the skein relations  (\ref{skeinsu2over}) and (\ref{skeinsu2under}).}

\begin{equation}\labell{skeinsu2over1}
\langle \overcrossing \rangle = -\Big(\cosh(\sqrt{3} \beta)\Big) \Big( (\smoothing) + (\hsmoothing)\Big)  
+ \Big(\frac{1}{\sqrt{3}}\sinh(\sqrt{3} \beta)\Big) \Big( -(\smoothing) + (\hsmoothing)\Big)  
\end{equation}

and

\begin{equation}\labell{skeinsu2under1}
\langle \undercrossing \rangle =  -\Big(\cosh(\sqrt{3} \beta)\Big) \Big( (\smoothing) + (\hsmoothing)\Big)  
- \Big(\frac{1}{\sqrt{3}}\sinh(\sqrt{3} \beta)\Big) \Big( -(\smoothing) + (\hsmoothing)\Big)  
.\end{equation}

\noindent Alternatively,

\begin{equation}\labell{skeinsu2over2}
\langle \overcrossing \rangle =  a (\smoothing) + b (\hsmoothing)
\end{equation}

\noindent and

\begin{equation}\labell{skeinsu2under2}
\langle \undercrossing \rangle = b (\smoothing) + a (\hsmoothing)
\end{equation}

\noindent where 

\begin{equation}\labell{a}
a = -\cosh(\sqrt{3} \beta) - \frac{1}{\sqrt{3}}\sinh(\sqrt{3} \beta)
\end{equation}

\noindent and

\begin{equation}\labell{b}
b = -\cosh(\sqrt{3} \beta) + \frac{1}{\sqrt{3}}\sinh(\sqrt{3} \beta).
\end{equation}

If we further normalize the expectations by a factor of 

$$(\sqrt{ab})^{\# {\rm crossings}}$$

\noindent the normalized expectation satisfies the relations\footnote{Again, we abuse notation by using the same notation for the expectation with the modified normalization.}

\begin{equation}\labell{k1}
\langle\overcrossing\rangle= q (\smoothing) + q^{-1} (\hsmoothing)\end{equation}
\begin{equation}\labell{k2}\langle\undercrossing\rangle  = q^{-1} (\smoothing) + q (\hsmoothing)\end{equation}
\noindent where 

$$q = \frac{a}{\sqrt{ab}}.$$

Note that a rotation by $\frac\pi{2}$ exchanges the right hand sides of (\ref{k1}) and (\ref{k2}); we are therefore justified in writing these equations in the form

\begin{equation}\labell{ku1}
\langle\slashoverback\rangle= q (\smoothing) + q^{-1} (\hsmoothing)\end{equation}
\begin{equation}\labell{ku2}\langle\backoverslash\rangle  = q^{-1} (\smoothing) + q (\hsmoothing)\end{equation}

\noindent which involve unoriented links.\footnote{Equations (\ref{ku1}) and (\ref{ku2}) are a type of deformation of the classical relation (\ref{su2miracle}), which is a linear algebraic fact about $SU(2).$}  

The relations (\ref{ku1}) and (\ref{ku2}) (which are really the same relation, rotated by $\pi/2$) are precisely those of Kauffman's bracket polynomial \cite{K2}.  By adding
(\ref{ku1}) and (\ref{ku2}) we obtain the relation

\begin{equation}\labell{kbracket}
\langle \slashoverback \rangle + \langle \backoverslash \rangle = (q + q^{-1}) \Big((\smoothing) + (\hsmoothing) \Big)
\end{equation}

\noindent of the $L$-polynomial, also due to Kauffman, which is a close relative of the Kauffman bracket polynomial.\footnote{Subtraction of  (\ref{ku2}) from (\ref{ku1}) gives the relation of Kauffman's Dubrovnik invariant; see \cite{K2}, Section VII.}

So far all our expressions, culminating with (\ref{ku1}) and (\ref{kbracket}), are relations between functions of the background gauge field $A.$
Now Kauffman in \cite{K2} shows that his bracket polynomial is a regular isotopy invariant of link diagrams, which is close relative of the Jones polynomial, if, in addition to the relation $(\ref{ku1}),$ we impose

\begin{equation}\labell{zpe}
\langle \bigcirc \rangle =- (q^2 + q^{-2}).
\end{equation}

This choice is necessary for the bracket polynomial to be invariant under the second Reidemeister move.\footnote{See \cite{K2}, p. 431, second equation.} This choice also determines, using (\ref{ku1}), the behavior of the expectations under the first Reidemeister move.  A further adjustment of the bracket polynomial, using an orientation of the link components, is necessary to obtain an ambient isotopy invariant related to the Jones polynomial of oriented links. Note that invariance under the third Reidemeister move is immediate due to the definition of $\langle \cdot \rangle$ in terms of expectations by the formula (\ref{exp1}).

Now the choice (\ref{zpe}) is not the natural one from the point of view of Chern-Simons Gauge Theory, where from the Lagrangian formulation, we would expect the expectation $\langle W_{C_1}\dots W_{C_n}\rangle$ to be taken at $A = 0,$  so that $\langle \bigcirc \rangle = -2.$\footnote{The sign is again due to our adjustment of the normalization.}   Thus, as it stands, the expectations of products of Wilson loops in Chern-Simons Gauge Theory in axial gauge are closely related to the Jones polynomial, but not identical with it.\footnote{Note that the only values of $q$ for which  $(q^2 + q^{-2})= 2$ holds are $q = \pm 1.$}  One interpretation of this fact is that the gauge anomalies of Chern-Simons Gauge Theory cause what may be seen as a shift in a quantity $\langle \bigcirc \rangle,$ in some ways resembling a shift in the zero point energy, of the theory in axial gauge versus that of covariant gauge.

This issue requires further investigation.  One or two comments, however:  One possible remedy for this ailment\footnote{Note that renormalizing the function  $W_C(\cdot)$ by multiplication by some fixed constant will spoil the skein relations, so this type of  renormalization would not appear to be useful.} is the inclusion (which we have already made) of a background gauge field $J.$  Repeated applications of the relation (\ref{ku1}) to any link diagram will result in a function of the background field given by a sum of terms each of which is the product of holonomies of this field about a number of unlinked, unknotted circles $c_1,\dots, c_l.$   We may now define

\begin{equation}\labell{jexp}
\langle W_{c_1}(J)\dots W_{c_l}(J)\rangle_{\cJ} = (-(q^2 + q^{-2}))^l\end{equation}

\noindent to obtain the desired invariant.  \footnote{As we noted, since the unlinked, unknotted circles $c_1,\dots,c_n$ actually involve jumps, it is convenient to choose the background field $J$ in the final evaluation to be the pullback of a connection from the plane.  This does not affect our ability to consider a definition like (\ref{jexp}).}

This has the defect that $\langle \cdot \rangle_{\cJ}$ is not defined for expectations of anything other than products of Wilson loops corresponding to unlinked unknotted circles, but this can be remedied  by defining the expectation to be $(- (q^2 + q^{-2}))^l$ where $l$ is the number of components, for {\em any} link. Another (not quite circular) remedy to this problem is to define $\langle \cdot \rangle_{\cJ}$ for all links by taking it to be the formal Chern-Simons functional integral in the variable $J,$ and then defining the invariant to be $$\langle\langle W_{C_1}(A+J)\dots W_{C_n}(A+J)\rangle_{A=0}\rangle_{\cJ}.$$ This addition of a background gauge field $A \to A + J,$ where both $A$ and $J$ are fields with Chern-Simons Lagrangian, is in some ways reminiscent of the work of Beasley-Witten \cite{bw}.  
 
\section{$GL(n)$ and $U(n)$}\labell{glnsec}

We now perform similar calculations for higher rank Lie groups.  We note first that the calculations of expectations in (\ref{linear}) and (\ref{quad}) involve only the Lie algebra $\fg,$ not the group $G,$ and so can be generalized at once to noncompact Lie groups.  We study the case of $G = GL(n,\C)$ with the trace for the Wilson loops taken in the defining representation, with Lie algebra $\fg = {\mathfrak gl}(n,\C)= M_n(\C),$ the Lie algebra of all complex $n \times n$ matrices.  Given a basis $\{e_\alpha \}$ for $\fg$ we continue to write $g_{\alpha\beta}= {\rm Tr~} e_\alpha e_\beta,$ where ${\rm Tr}$ is also the trace in the defining representation, and extend formulas (\ref{exp1}) and (\ref{op}) in this way to define the expectations of products of Wilson loops in the case of noncompact groups.\footnote{This despite the fact that the metric on $\fg$ is given by $g(X,Y) = \tr (X^*Y)$ for $X,Y \in M_n(\C).$  Our definition of $g_{\alpha\beta}$ corresponds to taking the trace, not the metric, in the action (\ref{axgaction}), and all subsequent formulas.  For the Lie algebra of a compact simple Lie group, the metric differs from the bilinear form given by the trace only by a sign; for a noncompact group, the two bilinear forms may not be so simply related.}

We choose as a $\C$-basis for $ {\mathfrak gl}(n,\C)$ the matrices $E_{ij}$ given by

$$(E_{ij})_{kl} = \delta_{ik}\delta_{jl}.$$

\noindent for $i,j \in \{1,\dots,n\}.$  To compute the expectation $\langle W_{C_1} \dots W_{C_n}\rangle,$ we again find formulas for the action of the operators $D_i$ on a product of Wilson loops, as in (\ref{exp1}), by computing the right hand sides
of (\ref{linear}) and (\ref{quad}) in the explicit basis for $\fg$ given by the $E_{ij}.$  We have, from (\ref{linear})
 
\begin{equation}\begin{aligned}\labell{lineargl}
 4 \sum_{\alpha,\beta}  (g^{-1})_{\alpha\beta} u(t^+ - t^-) \frac{\delta}{\delta A_\alpha^1(x,t^+)}\frac{\delta}{\delta A_\beta^2(x ,t ^-)}  W_C(A)&\\ = 2 \sum_{i,j} {\rm tr~} \hol_X(A) E_{ij} \hol_Y(A) E_{ji}= 2 \tr U^{12} \tr V^{21}
 \end{aligned}\end{equation}

\noindent for an overcrossing, where $U^{12}$ and $V^{21}$ are, as in Section \ref{su2sec}, the holonomies of the connection $A$ along the paths $X,Y$ obtained by removing the points $(x,t^\pm)$  from $C.$  Likewise in the quadratic case, we have (again for an overcrossing)
 
 \begin{equation}\begin{aligned}\labell{quadgl}
4 \sum_{\alpha,\beta}   u(t^+ - t^-) (g^{-1})_{\alpha\beta} \frac{\delta}{\delta A_\alpha^1(x,t^+)}\frac{\delta}{\delta A_\beta^2(x ,t ^-)}  \tr(U^{11}) \tr(V^{22})  &\\=2 \sum_{i,j} \tr(U^{11} E_{ij}) \tr(V^{22} E_{ji}) = 2 \tr(U^{11}V^{22}) 
 \end{aligned}\end{equation}

 \noindent where $U^{11}, V^{22}$ are as before 
 
 $$U^{11} = {\rm tr~}  \hol_{C,(x,t^+)}(A)$$
 
 $$V^{22} = {\rm tr~}  \hol_{C^\prime,(x,t^-)}(A).$$
 
 We also compute (again for an overcrossing)
 
 \begin{equation}\begin{aligned}\labell{lingl2}
 4\sum_{\alpha,\beta}  u(t^+ - t^-) (g^{-1})_{\alpha\beta}  \frac{\delta}{\delta A_\alpha^1(x,t^+)}\frac{\delta}{\delta A_\beta^2(x ,t ^-)}  \tr(U^{12}) \tr(V^{21}) =&\\
\sum_{i,j} u(t^+ - t^-) \Big (  (\tr E_{ij} U^{12} E_{ji} )( \tr V^{21}) 
 &\\+   (\tr  U^{12})( \tr  E_{ij}V^{21}E_{ji})  
 &\\+ 2  (\tr E_{ij} U^{12} )( {\rm tr~} E_{ji}  V^{21} )\Big) &\\=\frac12\Big( 2n \tr U^{12} \tr V^{21} + 2 \tr U^{12} V^{21}\Big)
 \end{aligned}\end{equation}

\noindent and likewise
 
 \begin{equation}\begin{aligned}\labell{quadgl2}
4 \sum_{\alpha,\beta}  u(t^+ - t^-) (g^{-1})_{\alpha\beta} \frac{\delta}{\delta A_\alpha^1(x,t^+)}\frac{\delta}{\delta A_\beta^2(x ,t ^-)}  \tr(U^{11}V^{22})  &\\=\frac{1}{2} (2 \tr U^{11}\tr V^{22} + 2n \tr U^{11} V^{22}). 
   \end{aligned}\end{equation}
   
 So the action of the operator
 
 $$ D= 4 \sum_{\alpha,\beta}  u(t^+ - t^-) (g^{-1})_{\alpha\beta} \frac{\delta}{\delta A_\alpha^1(x,t^+)}\frac{\delta}{\delta A_\beta^2(x ,t ^-)}$$
 
 \noindent in either case is given, for an overcrossing, by
 
 \begin{equation}
 D
  \begin{bmatrix}
      \overcrossing    \\
   \upupsmoothing\ 
  \end{bmatrix}
 =
\tilde{M}_n
   \begin{bmatrix}
     \overcrossing  \\
  \upupsmoothing 
  \end{bmatrix}
\end{equation}

\noindent and for an undercrossing by

\begin{equation}
 D
  \begin{bmatrix}
   \undercrossing    \\
  \upupsmoothing 
  \end{bmatrix}
 =
-\tilde{M}_n
   \begin{bmatrix}
   \undercrossing   \\
  \upupsmoothing  
  \end{bmatrix}
\end{equation}

\noindent where $\tilde{M}_n$ is the $2 \times 2$ matrix 

\begin{equation}\labell{tmn}
\tilde{M}_n = \begin{bmatrix}
    0 & 1   \\
   2 & n
  \end{bmatrix}\end{equation}
  
The traceless part of $\tilde{M}_n$ is       

$$M_n = \tilde{M}_n -\frac{n}2  I = \begin{bmatrix}
   -\frac{n}2 & 1   \\
   2 & \frac{n}2
  \end{bmatrix}$$
  
  \noindent so that $M_2  = M,$ and our results agree with those we obtained for $SU(2)$ in Section \ref{su2sec}.  (See Section \ref{unsec} for more discussion of this point.)\footnote{Note that for $n = 1,$ we have$ \langle \undercrossing \rangle = \langle\overcrossing\rangle = \langle \upupsmoothing \rangle,$ so that $M_1$ can be considered as the $1 \times 1$ matrix $M_1 = 2,$ in agreement with the results of Section \ref{u1sec}.}
  
The matrix $M_n$ is again a $2 \times 2$ traceless matrix, so that

$$M_n^2 = |\det(M_n)| I = (\frac{n^2}4 + 2)I .$$

To evaluate expectations $\langle W_{C_1}\dots W_{C_n}\rangle,$ we again have to exponentiate $\tilde{M}_n,$ and we have\footnote{See also \cite{helgason}, p. 149, Exercise B1.}
\begin{equation} \labell{expmn}
e^{\beta\tilde{M}_n} = e^{\beta \frac{n}2 } e^{\beta M_n} = 
e^{\beta \frac{n}2 }\Big( \cosh(\beta {\sqrt{\Delta_n}} ) I + \frac{1}{\sqrt{\Delta_n}} \sinh(\beta{\sqrt{\Delta_n}} )M_n\Big)
 \end{equation}
 
 \noindent where $\Delta_n =  \frac{n^2}4 + 2.$
 
Thus
 
 \begin{equation}\labell{oskeingln} 
 e^{-\beta \frac{n}2 }\langle \overcrossing \rangle = \Big( \cosh(\beta {\sqrt{\Delta_n}} ) -\frac{n}2  \frac{1}{\sqrt{\Delta_n}} \sinh(\beta{\sqrt{\Delta_n}} )\Big) (\doublepoint)  + \frac{2}{\sqrt{\Delta_n}} \sinh(\beta{\sqrt{\Delta_n}} ) (\upupsmoothing)
 \end{equation}
 
 \noindent while
 
  \begin{equation}\labell{uskeingln} 
 e^{\beta \frac{n}2 }\langle \undercrossing \rangle = \Big( \cosh(\beta {\sqrt{\Delta_n}} ) + \frac{n}2  \frac{1}{\sqrt{\Delta_n}} \sinh(\beta{\sqrt{\Delta_n}} )\Big) (\doublepoint)  - \frac{2}{\sqrt{\Delta_n}} \sinh(\beta{\sqrt{\Delta_n}} ) (\upupsmoothing)
 \end{equation}
 
 Note that the framing terms 
 
 $$ e^{\pm \frac{n}2 \beta}$$
 
 \noindent arise from the trace of $\tilde{M}_n.$  We would expect that the $SL(n,\C)$ expectations differ from the $GL(n,\C)$ expectations by these terms.
 
 We now use equations  (\ref{oskeingln}) and (\ref{uskeingln}) to obtain skein relations in the form usual for the HOMFLY invariants \cite{homfly}.  Subtracting  (\ref{uskeingln}) from (\ref{oskeingln}) we have
 
  \begin{equation}\labell{-skeingln} 
 e^{-\beta \frac{n}2 }\langle \overcrossing \rangle - e^{\beta \frac{n}2 }\langle \undercrossing \rangle= \Big(- \frac{n}{\sqrt{\Delta_n}} \sinh(\beta{\sqrt{\Delta_n}} )\Big) (\doublepoint)  + \frac{4}{\sqrt{\Delta_n}} \sinh(\beta{\sqrt{\Delta_n}} ) (\upupsmoothing)
 \end{equation}
 
 \noindent while adding equations  (\ref{oskeingln}) and (\ref{uskeingln}) we obtain

  \begin{equation}\labell{+skeingln} 
 e^{-\beta \frac{n}2 }\langle \overcrossing \rangle + e^{\beta \frac{n}2 }\langle \undercrossing \rangle= 
 \Big(2 \cosh(\beta {\sqrt{\Delta_n}} )  \Big) (\doublepoint)  
 \end{equation}
 
We may use (\ref{+skeingln}) to eliminate the $\doublepoint$ term from (\ref{-skeingln})  and obtain the HOMFLY type skein relation

\begin{equation}\labell{hskeingln} 
 \Big( 1+ \frac{n}{2\sqrt{\Delta_n}} \tanh (\beta \sqrt{\Delta_n})\Big)
 e^{-\beta \frac{n}2 }\langle \overcrossing \rangle -\Big( 1- \frac{n}{2\sqrt{\Delta_n}} \tanh (\beta \sqrt{\Delta_n})\Big)
 e^{\beta \frac{n}2 }\langle \undercrossing \rangle= 
 \frac{4}{\sqrt{\Delta_n}} \sinh(\beta{\sqrt{\Delta_n}} ) (\upupsmoothing)
 \end{equation}
 
 It is of course possible to divide (\ref{hskeingln}) by 
 
 $$\sqrt{1- \frac{n^2}{4\Delta_n} \tanh^2(\beta \sqrt{\Delta_n})}$$
 
 \noindent to obtain a skein relation of a more standard type
 
 $$q \langle \overcrossing \rangle - q^{-1}\langle \undercrossing \rangle = z ( \upupsmoothing ).$$
 
 We note that although \cite{homfly} guarantees that there are HOMFLY link invariants satisfying the skein relation (\ref{hskeingln}), this does not imply that our expectations, although they also satisfy this relation, are equal to this HOMFLY invariant.  We expect instead that adjustments of various sorts similar to those of Section \ref{su2sec} will be needed to relate the two.  In particular, the expectations satisfy both (\ref{+skeingln}) and (\ref{-skeingln}), not just the one relation (\ref{hskeingln}).  We hope to return to these issues in a future paper.
 
 \subsection{$G = U(n)$}\labell{unsec}
 
Consider now again the $GL(n,\C)$ case, but replace the basis $E_{ij}$ with the basis
given by the matrices

$$ \frac{E_{ij}- E_{ji}}{\sqrt{2}} , \sqrt{-1} \frac{E_{ij} + E_{ji}}{\sqrt{2}}$$

\noindent  where $i > j,$ alongside the matrices $\sqrt{-1} E_{ii}.$   These matrices form a $\C$-basis for the Lie algebra of $GL(n,\C),$ and therefore give the same relations for expectations of products of Wilson loops in the $GL(n,\C)$ case derived previously using the basis $E_{ij}.$\footnote{This can also be seen by explicit computation.}  On the other hand, these matrices also give an $\R$-basis for the Lie algebra of $U(n).$  So applying formulas (\ref{linear}) and (\ref{quad}), and using this basis, we obtain the same formulas in the $U(n)$ case (in the defining representation) as in the $GL(n,\C)$ case.  It follows therefore that the $U(n)$ expectations satisfy the same skein relations (\ref{-skeingln}) and (\ref{+skeingln}) as the $GL(n,\C)$ expectations.  This also explains the relation $M = M_2$ we noted above.

\section{The Goldman Bracket}\labell{goldmansec}

We now perform a computation in two dimensional gauge theory similar to those we have been doing with the symmetric operator given by the propagator $\Gamma$ on the space of connections on the three-manifold $\R^3.$  If we consider the space of connections on a two-manifold $\Sigma,$ the Atiyah-Bott symplectic form gives a Poisson bivector field $\Pi,$ that is, a skew operator, whose action on functions given by Wilson loops can be computed by methods similar to ours.  We will see that these computations give results identical to those given by Goldman \cite{goldman} on the moduli space of flat connections on the two manifold.  Since this moduli space is the reduced space of the space of connections by the action of the gauge group, and since the Wilson loops are invariant under this group, this is to be expected.

In a companion paper \cite{star} we will show that the similarity between Goldman's formulas and the link expectation formulas in this paper are part of a larger relation between deformation quantization of spaces of connections and Chern-Simons gauge theory, echoing the relation between Chern-Simons gauge theory and geometric quantization of moduli spaces of flat connections on a two manifold.

So let $\Sigma$ be a closed, oriented two manifold, let $G$ be a Lie group as before, along with a choice of representation whose associated trace gives a metric on $\fg,$ and let $\cA(\Sigma)$ be the space of connections on the trivialized principal $G$-bundle on $\Sigma.$  Then $\cA(\Sigma) = \Omega^{1}(\Sigma)\otimes \fg.$  This space is equipped \cite{ab} with a symplectic form $\Omega$ given by

$$ \Omega_A (v,w) =  {\rm Tr} \int_\Sigma v \wedge w$$

\noindent  where $A \in \cA(\Sigma)$ and $v,w \in T\cA_A = \Omega^{1}(\Sigma)\otimes \fg,$ and the trace ${\rm Tr}$ is the trace in the chosen representation.

The Poisson bivector field $\Pi$ arising from the symplectic form $\Omega$ is given in local coordinates  by 

\begin{equation}\labell{poiss}\Pi =   \int_\Sigma dx \sum_{\alpha,\beta} \sum_{a,b} \epsilon_{ab}   (g^{-1})_{\alpha\beta} \frac{\delta}{\delta A_\alpha^a(x)}\frac{\delta}{\delta A_\beta^b(x)}\end{equation}
 
 \noindent which is a close relative of the operators $D_i$ in (\ref{op}).  Here again we choose a basis $\{e_\alpha\}$ for $\fg,$ and write $g_{\alpha\beta}= {\rm Tr~} (e_\alpha e_\beta)$ and local coordinates on $\Sigma$ near each point $x,$ in terms of which we express the components of $A$ as a section of the cotangent bundle.\footnote{An invariant expression for $\Pi$ is given as follows.  Given a differentiable function $f$ on $\cA(\Sigma),$ the derivative is $\delta f |_A \in \Omega^1(\Sigma,\fg)^*.$  The orientation of $\Sigma$ gives a map $\Phi: \Omega^1(\Sigma,\fg)\to \Omega^1(\Sigma,\fg)^*$ given by $\Phi(v)(w) = {\rm Tr} \int_\Sigma v \wedge w$ for $v,w \in \Omega^1(\Sigma,\fg).$  In these terms, morally $\Pi (fg) (A) ={\rm Tr} \int_\Sigma \Phi^{-1}(\delta f|_A) \wedge  \Phi^{-1}(\delta g|_A),$ for any differentiable functions $f, g$ on $\cA(\Sigma)$ whose derivatives are in the image of $\Phi.$  (In general, these derivatives may correspond under $\Phi$ to distributional one forms.) This coordinate free definition agrees with (\ref{poiss}).}
  
Suppose that we are given two oriented simple closed curves $C,C^\prime$ in $\Sigma,$ intersecting transversally at a finite number of points $x_1,\dots,x_n.$   We wish to compute the Poisson bracket

$$\{ W_C, W_{C^\prime}\}(A)$$

\noindent (where as before we write $W_C(A) = {\rm tr}_R \hol_{C,\star}(A)={\rm tr~} \hol_{C,\star}(A)$ for the trace in some representation $R$ of the holonomy of the connection $A$ along the oriented path given by $C$ with basepoint $\star,$ and suppress for notational simplicity the basepoint where we are taking a trace, as well as the representation) by applying the skew operator $\Pi$ to the pair of functions $W_C,W_{C^\prime}.$

We obtain\footnote{Strictly speaking, we would need to smooth the functions $W_C.$ To do this, smooth the connection $A$ by convolving with a delta sequence $\delta_\epsilon,$ and write $W_C^\epsilon(A) = W_C(A \star \delta_\epsilon).$ Applying $\Pi$ to the product $W_C^\epsilon W_{C^\prime}^\epsilon $ of two such functions and taking the limit $\epsilon \to 0,$  we obtain (\ref{pbi}).}

\begin{equation}\labell{pb}\{ W_C , W_{C^\prime} \} (A)=  \sum_i \sum_{\alpha,\beta}  \sum_{a,b} (g^{-1})_{\alpha\beta}  \epsilon_{ab}  \frac{\delta W_C(A)}{\delta A_\alpha^a(x_i)}\frac{\delta W_{C^\prime}(A)}{\delta A_\beta^b(x_i)}  .
\end{equation}

\noindent where the $x_i$ are the transverse intersection points of $C$ and $C^\prime.$  Then, as in the three dimensional case (\ref{quad}), the contribution of the $i$-th point to this sum is

\begin{equation}\labell{pbi}\sum_{\alpha,\beta} (g^{-1})_{\alpha\beta}  \epsilon_i  \tr( \hol_{C,x_i}(A) e_\alpha) \tr( \hol_{C^\prime,x_i}(A) e_\beta)\end{equation}

\noindent where, as before, $e_\alpha$ is a basis for $\fg$ and $\epsilon_i = \pm 1$ is a sign arising from the orientation on $\Sigma$ and the orientation of the curves $C,C^\prime$ ; explicitly, if $v_i,v_i^\prime \in T\Sigma|_{x_i}$ are tangents to the oriented curves $C,C^\prime$ at $x_i,$ pointing in the direction of their orientations, then\footnote{Note that if the curves $C,C^\prime$ do not intersect, the corresponding Wilson loops Poisson commute; see \cite{weitsman} for a proof of this fact along the lines of the gauge theoretic computation given here.}

$$\epsilon_i = {\rm sgn}( \omega_{x_i}(v_i,v_i^\prime)),$$

\noindent where $\omega$ is a symplectic form on $\Sigma$ with $\int_\Sigma \omega > 0.$

If we choose the same representation for the Wilson loops as we did in defining the symplectic form, so that ${\rm tr} = {\rm Tr},$ we have

\begin{equation}\labell{gt} \sum_{\alpha,\beta} (g^{-1})_{\alpha\beta}  \tr(U e_\alpha) \tr (V e_\beta) = \tr (\pi(U) \pi(V))\end{equation} 

\noindent where $\pi : G \to \fg$ is the map given by composing the inclusion of $G$ into $GL(m,\C)$ given by the chosen representation, followed by the restriction to $GL(m,\C) \subset {\mathfrak gl}(m,\C)$ of the projection ${\mathfrak gl}(m,\C) \to \fg$ in the metric arising from the trace ${\rm Tr}.$

In the case where $G = GL(n,\C)$ and the representation is the defining representation, the map $\pi$ is the identity, and we obtain\footnote{This may also be checked by direct computation using the basis $\{E_{ij}\}$ from Section \ref{glnsec}.}

\begin{equation}\labell{gtgln}\sum_{\alpha,\beta} (g^{-1})_{\alpha\beta}  \tr(U e_\alpha) \tr (V e_\beta)= \tr (UV).\end{equation}

Applying this to the Poisson bracket $\{ W_C , W_{C^\prime} \},$ equations (\ref{pb}) (\ref{pbi}) and (\ref{gtgln}) give

\begin{equation}\labell{pbgln}\{ W_C , W_{C^\prime} \} (A)= \sum_i \epsilon_i W_{C \star_i C^{\prime}}(A)\end{equation}

\noindent where ${C \star_i C^{\prime}}$ denotes the concatenation of the oriented curves $C$ and $C^{\prime}$ at the common basepoint $x_i.$  This agrees with Goldman's computation (Theorem 3.13 of \cite{goldman}).

In the case of $G = SU(2)$ and the defining representation, the map $\pi:SU(2) \to su(2)$ is given
by 

$$\pi(U) = \frac12(U - U^{-1}).$$

\noindent (cf. Corollary 1.10 of \cite{goldman}).  We therefore have

\begin{equation}\labell{pbsu2}\{ W_C , W_{C^\prime} \} (A)= \frac12 \sum_i \epsilon_i (W_{C \star_i C^{\prime}}(A) - W_{C \star_i \bar{C}^{\prime}}(A) )\end{equation}

\noindent  where we denote by $\bar{C}$ the curve $C$ with the orientation reversed.  This is in agreement with Theorem 3.17 of \cite{goldman}.\footnote{ The Goldman bracket was the origin of the subject of String Topology \cite{cs}, which generalized Goldman's construction on two manifolds to homology of free loop spaces.  The operator $\Gamma$ in (\ref{propagator}) is a type of even analog of the Poisson bivector field $\Pi,$ and the skein relations are in this sense the analogs of the products in Goldman's Lie algebra of curves.  It would be interesting to investigate if there is an even analog of String Topology bearing the same relation to $\Gamma$ as String Topology does to $\Pi.$}

\end{document}